\documentclass[preprint,12pt]{elsarticle}

\usepackage{graphicx}
\usepackage{epsfig}

\usepackage{amssymb}

\journal{Linear Algebra and its Applications}

\begin{document}

\begin{frontmatter}

\title{A Lower Bound for Algebraic Connectivity based on Connection Graph Stability Method}


\author[]{Ali Ajdari Rad\footnotemark[1], Mahdi Jalili, and Martin Hasler}
\address{Laboratory of Nonlinear Systems, \'{E}cole Polytechnique F\'{e}d\'{e}ral de Lausanne, School of
Computer and Communication Sciences, 1015 Lausanne, Switzerland}
\address{\{ali.ajdarirad, mahdi.jalili, martin.hasler\}@epfl.ch}

\begin{abstract}
In this paper a tight lower bound for algebraic connectivity of
graphs (second smallest eigenvalue of the Laplacian matrix of the
graph) based on \emph{connection-graph-stability} method is introduced. The connection-graph-stability score for each edge is defined as the sum of the length of
all the shortest paths making use of that edge. We prove that the algebraic connectivity of the graph is lower bounded by the size of the graph divided by the maximum connection graph stability of the edges.
\end{abstract}

\begin{keyword}
Algebraic connectivity \sep Graph Laplacian \sep Connection graph stability score
% PACS codes here, in the form: \PACS code \sep code PATREC-D-04-00385
%\PACS
\end{keyword}
\end{frontmatter}

% main text
\newenvironment{definition}[1][Definition]{\begin{trivlist}
\item[\hskip \labelsep {\bfseries #1}]}{\end{trivlist}}

\section{Introduction}
\label{sect:intro}  % \label{} allows reference to this section
\footnotetext[1]{Corresponding author.}
  
Let $G = (V,E)$ be a connected simple graph with $n = |V|$
vertices and $|E|$ edges. For a graph $G$ a simple path with minimum length connecting two
nodes $u$ and $v$ is called shortest path between the nodes and
denoted by $P_{uv}$. The longest shortest path between pairs of
nodes is called \emph{diameter} of the graph, denoted by
$D_{max}$. The Laplacian matrix of $G$ is defined as
$L = D - A $, where $A$ is the binary adjacency matrix and
$D=diag(d_u; u \in V )$ is the degree-diagonal matrix of $G$. $L$ is a
positive semidefinite, symmetric and singular matrix whose
eigenvalues are in the form of
$\lambda_{n}(G)\geq\lambda_{n-1}(G)\geq\ldots\lambda_{2}(G)\geq\lambda_{1}(G)=0$.
These eigenvalues are important in graph theory and have close relations to numerous graph invariants. Among them, $\lambda_2$, also called algebraic connectivity, has attracted more attention. There were some attempts to find a lower bound for $\lambda_2$ based on the simple properties of graphs such as diameter, order, and number of edges of the graph (see \cite{Abreu07} for a comprehensive review). For example, Mohar \cite{Mohar91} showed that $\lambda_2\geq\frac{4}{nD_{max}}$ %Alon and Milman [1] obtained $D_{max}\leq[\sqrt{2k_{max}/\lambda_2}\log{n}_2]$,
and recently Lu \cite{Lu07} proved that $\lambda_2\geq \frac{2n}{2 + (n-1)nD_{max}-2|E|D_{max}}$. In this paper, we present a new lower bound for $\lambda_2$ based on the connection-graph-stability scores associated to the edges, which are defined for each edge as the sum of the length of all the shortest paths making use of that edge. Then, we will prove that the proposed lower bound is always tighter than the previously mentioned bounds proposed by Mohar and Lu.

\section{Lower Bound based on Connection-graph-stability method}
The connection-graph-stability method was proposed by Belykh et al \cite{Belykh04a} to establish a criteria for the global stability of synchronization manifold of a network of coupled dynamical systems. Here, we will reuse the concept to obtain a lower bound for algebraic connectivity of a graph.
\newtheorem{thm}{Theorem}
\newtheorem{theorem3}{Definition}
\begin{theorem3}[Connection-graph-stability]
If for each pair of nodes in the graph a simple connection path is
considered (not necessarily the shortest path), a Connection-graph-stability score for each edge $k$ of the graph is denoted by $C_k$
and defined as the sum of the length of the paths passing through the edge, i.e.
 \begin{equation}\label{frm0}
C_k = \frac{1}{2}\sum_{u=1}^{n}\sum_{v=1}^{n}\varphi_{uv}(k)|P_{uv}|,
\end{equation}

where

\end{theorem3}

\begin{displaymath}
\varphi_{uv}(k)=\left\{\begin{array}{lcrr}
1 & & if & k \in P_{uv}\\
0 &  & & otherwise
\end{array}
\right.
\end{displaymath}

\begin{thm}
 Let $G$ be a connected simple graph with $n$ nodes. Then,
 \begin{equation}\label{frm1}
\lambda_{2}\geq \frac{n}{C_{max}},
\end{equation}
where $C_{max}$ is the maximum connection-graph-stability score
assigned to the edges of $G$, i.e. $C_{max} = \max_{k \in E}{C_k}$.
\end{thm}
 \textbf{Proof.} Let $G$ be a simple graph. Fiedler \cite{Fiedler75} showed that

\begin{equation}\label{frm2}\lambda_2=\min \frac{2n \sum_{uv\in
E}(x_{u}-x_{v})^2}{\sum_{u\in V}\sum_{v\in
V}(x_{u}-x_{v})^2},\end{equation}

where the minimum is taken over all non-constant vectors $x=(x_{v})_{v\in V(G)}$ with $\|x\|=1$.\\
Denoting $X_{uv}=x_{u}-x_{v}$, equation (\ref{frm2}) can be written as\\

 \begin{equation}\label{frm3}
    \lambda_2=\min \frac{2n \sum_{uv\in
E}(X_{uv})^2}{\sum_{u\in V}\sum_{v\in V}(X_{uv})^2}.
\end{equation}

Let $P_{uv}=um_1m_2...m_kv$ with $m_i\in V$ being the path connecting $u$ to $v$ and $|P_{uv}|$ its length. $X_{uv}$ can be
expressed as \[X_{uv} = X_{um_1}+X_{m_1m_2}+...+X_{m_k
v}=\sum_{e\in P_{uv}}X_e.\] Applying the Cauchy-Schwartz inequality,
one obtains
\[X_{uv}^2 = (\sum_{e\in P_{uv}}1.X_e)^2\leq |P_{uv}|\sum_{e\in
    P_{uv}}X_e^2.\]

Then,

 \begin{equation}\label{frm4}\sum_{u=1}^n\sum_{v=1}^n(X_{uv})^2\leq
 \sum_{u=1}^n\sum_{v=1}^{n}(|P_{uv}|\sum_{e\in E} \varphi_{uv}(e)X_e^2)=\sum_{e\in E}2C_eX_e^2\leq 2C_{max}\sum_{e\in E}X_e.\end{equation}

Substituting (\ref{frm4}) in (\ref{frm3})
\[\lambda_2=\min \frac{2n \sum_{uv\in
E}(x_{u}-x_{v})^2}{\sum_{u\in V}\sum_{v\in V}(x_{u}-x_{v})^2}\geq
\min \frac{2n \sum_{uv\in E}(x_{u}-x_{v})^2}{2C_{max}\sum_{uv\in
E}(x_{u}-x_{v})^2},\] and finally $\lambda_2 \geq \frac{n}{C_{max}}$. \\ \\ $\square$

\begin{theorem3}[Alternative Paths] The alternative
paths of any path $P_{uv}$ between two vertices $u$ and $v$ are the
paths connecting $u$ and $v$ with length less or equal to
$P_{uv}$, e.g. all shortest paths between $u$ and $v$ are
alternative paths for each other. The number of all possible
shortest paths between $u$ and $v$ is called the number of
alternative shortest paths between $u$ and $v$ and denoted by
$n_{uv}$. We denote the set of all possible shortest paths
between $u$ and $v$ by $\mathbb{P}_{uv}=\{P_{uv}^{(1)},
P_{uv}^{(2)}, ..., P_{uv}^{(n_{uv})}\}$
(\emph{i.e.}$n_{uv}=|\mathbb{P}_{uv}|)$.
\end{theorem3}

\begin{theorem3}[Path weighting strategy]
For any pair $u, v\in V$ choose a vector $\alpha_{uv}=(\alpha_{uv}^{(1)},\ldots, \alpha_{uv}^{(n_{uv})})$ such that $\alpha_{uv}^{(q)}\geq 0$ and $\sum_{q=1}^{n_{uv}}\alpha_{uv}^{(q)}=1$. The corresponding path weighting strategy is denoted by $\alpha$, the set of all these vectors.
\end{theorem3}

\begin{theorem3}[Extended connection-graph-stability score]
Then the extended connection graph stability score $C_k(\alpha)$ for edge $k$ is the sum of the weighted lengthes of all shortest paths making use of edge $k$, i.e.

 \begin{equation}\label{frm4.5}
C_k(\alpha) = \frac{1}{2}\sum_{u=1}^{n}\sum_{v=1}^{n}|P_{uv}|\sum_{q=1}^{n_{uv}}\varphi_{uv}^{(q)}(k)\alpha_{uv}^{(q)},
\end{equation}
where
\begin{displaymath}
\varphi_{uv}^{(q)}(k)=\left\{\begin{array}{lcrr}
1 & & if & k \in P_{uv}^{(q)}\\
0 & & & otherwise
\end{array}
\right.
\end{displaymath}
and $\alpha$ is the corresponding path weighting strategy.
\end{theorem3}

\begin{thm}
 Let $G$ be a connected simple graph with $n$ nodes. Then, for any path weighting strategy $\alpha$ we have
\begin{equation}\label{frm5}
\lambda_{2}=a(G)\geq \frac{n}{C_{max}(\alpha)},
\end{equation}
where $C_{max}(\alpha)$ is the maximum extended connection-graph-stability score assigned to the edges of $G$.
\end{thm}
 \textbf{Proof.} Consider $\mathbb{P}_{uv}$ and $n_{uv}$ for two arbitrary vertices.
Let $P_{uv}^{(q)} =um_1m_2...m_kv$ (with $m_i\in V(G)$) be the
$q_{th}$ path that connects $u$ to $v$ with length
$|P_{uv}^{(q)}|=|P_{uv}|$. Using the $q_{th}$ shortest path,
$X_{uv}$, can be expressed as

 $X_{uv} = x_u-x_v= X_{um_1}+X_{m_1m_2}+...+X_{m_kv}=\sum_{e\in
P_{uv}^{(q)}}X_e$.\\ \\
 Using the Cauchy-Schwartz inequality
  \begin{equation}\label{frm6}
(X_{uv})^2 = (\sum_{e\in P_{uv}^{(q)}}1.X_e)^2\leq
|P_{uv}|\sum_{e\in P_{uv}^{(q)}}(X_e)^2.
\end{equation}
On the other hand, one can express $X_{uv}$ as weighted average of
its alternative expansions as follows
\begin{equation}\label{frm7}
X_{uv}^2 = \sum_{q=1}^{n_{uv}}\alpha_{uv}^{(q)}
(X_{uv})^2\leq
\sum_{q=1}^{n_{uv}}\alpha_{uv}^{(q)}|P_{uv}|
\sum_{e\in P_{uv}^{(q)}}(X_e)^2,\end{equation}

Hence,
\[\sum_{u=1}^n\sum_{v=1}^n(X_{uv})^2\leq
 \sum_{u=1}^n\sum_{v=1}^{n}\sum_{q=1}^{n_{uv}}(\alpha_{uv}^{(q)}|P_{uv}|\sum_{e\in P_{uv}^{(q)}}X_e^2)\]
\[= \sum_{u=1}^n\sum_{v=1}^{n}|P_{uv}|\sum_{q=1}^{n_{uv}}\sum_{e\in E}(\varphi_{uv}^{(q)}(e)\alpha_{uv}^{(q)}X_e^2)
=\sum_{e\in E}2C_e(\alpha)X_e^2\]
\begin{equation}\label{frm9}
\leq 2C_{max}(\alpha)\sum_{e\in E}X_e^2.
\end{equation}

Substituting (\ref{frm9}) in (\ref{frm3}), we obtain
\[\lambda_2=\min \frac{2n \sum_{uv\in
E}(X_{uv})^2}{\sum_{u\in V}\sum_{v\in V}(X_{uv})^2}\geq \min
\frac{2n \sum_{uv\in E}(X_{uv})^2}{2C_{max}(\alpha)\sum_{uv\in
E}(X_{uv})^2},\] and finally $\lambda_2 \geq
\frac{n}{C_{max}(\alpha)}.$\\$\square$

\newtheorem{theorem4}{Corollary}
\begin{theorem4}
For any connected simple graph $G$ and any path weighting strategy $\alpha$ we have
\[\frac{n}{C_{max}(\alpha)}\geq \frac{n}{C_{max}}.\]
\end{theorem4}
\textbf{Proof.} Proof is straight forward since the connection-graph-stability score is a special case of the extended connection-graph-stability score.
 $\square$

\section{Comparing with the other lower bounds}
\subsection{Mohar's lower bound}
 Mohar \cite{Mohar91} showed that for any connected simple graph of order $n$ and diameter $D_{max}$
  \begin{equation}\label{Mohar}
    \lambda_2\geq \frac{4}{nD_{max}}.
\end{equation}

\begin{thm}
For any connected graph $G$ with $n$ vertices, diameter $D_{max}$,
and maximum connection-graph-stability number $C_{max}$, we have
\[\frac{n}{C_{max}}\geq \frac{4}{nD_{max}}.\]
\end{thm}
\textbf{Proof.} Consider the set of all shortest paths passing through
an edge $v_1v_2$. Define two subsets of vertices as follows: $v\in V_1$ if there is a
shortest path from $v_2$ to $v$ containing $v_1$ and $v\in V_2$ if there is a shortest path from $v_1$ to $v$ containing $v_2$. Note that $v_1\in V_1$ and $v_2\in V_2$. Then $V_1\cap V_2=\emptyset$.

If this is not the case, there is a vertex $v\in V_1\cap V_2$. Then there is shortest path $vP_1v_2v_1$ from $v$ to $v_1$ and a shortest path $vP_2v_1v_2$ from $v$ to $v_2$. If $|P_1|\leq|P_2|$ then the path $vP_1v_2$ is a shorter path from $v$ to $v_2$ than $vP_2v_1v_2$ which is a contradiction. If $|P_2|\leq|P_1|$ we come to the same contradiction.

Now let $m_1 = |V_1|$ and $m_2=|V2|$, then there are at most $m_1m_2$ shortest paths
between $V_1$ and $V_2$ passing through $e$, hence $C_e\leq
m_1m_2D_{max}$, where $C_e$ denotes the connection-graph-stability score of edge $e$. In addition, $m_2\leq n-m_1$, thus $C_e\leq m_1(n-m_1)D_{max}$, which is maximized for $m_1 = \frac{n}{2}$, therefore for any edge $e$ of the graph, and hence for the edge that has the maximum connection-graph-stability score, we have
\[C_e\leq (\frac{n}{2})^2D_{max}\Rightarrow \frac{n}{C_e} \leq \frac{n}{(\frac{n}{2})^2D_{max}}=\frac{4}{nD_{max}}.\]

$\square$

\subsection{Lu's lower bound}
 Lu \cite{Lu07} obtained the following bound for $\lambda_2$ of the Laplacian of any connected simple graph of order $n$, the number of edges $|E|$ and the diameter
 $D_{max}$
  \begin{equation}\label{Lu}
    \lambda_2\geq \frac{2n}{2 + (n-1)nD_{max}-2|E|D_{max}}.
\end{equation}

\begin{thm}
For any connected graph $G$ with $n$ vertices, number of edges
$|E|$, diameter $D_{max}$ and maximum connection-graph-stability
score $C_{max}$, we have
\[\frac{n}{C_{max}}\geq \frac{2n}{2 + (n-1)nD_{max}-2|E|D_{max}}.\]
\end{thm}
\textbf{Proof.} The total number of the shortest paths in $G$ is
$\frac{n(n-1)}{2}$ where out of these paths, there are $|E|$ ones
with length equal to one. Suppose that $e$ is the edge corresponding to $C_{max}$. There is only one path of length one passing
through this edge (the path that connects adjacent vertices of
$e$). At the same time, at most $\frac{n(n-1)}{2}-|E|$ paths
with length more than one can make use of $e$. According to the
definition of the $C_{max}$, the connection-graph-stability score of
$e$ is equal to the sum of the length of these paths plus the
length of the path connecting adjacent vertices of $e$, which is
one. Recall the maximum path length, i.e. diameter $D_{max}$, thus,
\[C_{max}\leq 1 + (\frac{n(n-1)}{2}-|E|)D_{max}\]
\[\Rightarrow \frac{n}{C_{max}}\geq \frac{2n}{2 + (n-1)nD_{max}-2|E|D_{max}}.\]
 $\square$

\section{Application to some well-known graphs} The maximum connection-graph-stability
score of some well-known graphs can be calculated analytically \cite{Belykh04a}, thus
the proposed bound can also be calculated for such graphs. Table 1 summarizes the
results on Complete, Path, Cycle, Star, and Peterson graphs.\\
\begin{table}[h]
\caption{Algebraic connectivity, maximum connection-graph-stability
number and proposed lower bound for some well-known graphs}
\begin{tabular}{lcccc}
  \hline
  Graph & $\lambda_2$ & Mohar & Lu & The Lower Bound (1)  \\
  \hline
  Complete graph & n & $\frac{4}{n}$ & n & n \\
  Path (n is even)& $2(1-\cos(\frac{\pi}{n}))$ & $\frac{4}{n(n-1)}$ & $\frac{2n}{2 + (n-2)(n-1)^2}$ & $\frac{8}{n^2}$  \\
  Cycle (n is odd)& $2(1-\cos(\frac{2\pi}{n}))$ & $\frac{8}{n(n-1)}$ & $\frac{2n}{2 + n(n-1)(\frac{n-3}{2})}$ & $\frac{24n}{(n^2-1)}$ \\
  Star  & 1 & $\frac{2}{n}$ &  $\frac{n}{1 + (n-2)(n-1)}$ & $\frac{n}{2n-3}$\\
  Peterson graph & 2 & 0.2 & 0.164 & 1.11\\
  \hline
\end{tabular}
\end{table}

\section{Conclusion}
In this paper a novel lower bound for algebraic connectivity of graphs based on the connection-graph-stability method was presented. It was proved that the maximum connection-graph-stability score which is defined as the maximum sum of the length of all the shortest paths making use of an edge is inversely related to the algebraic connectivity of the graph. In addition, it is proved that the proposed lower bound is always larger than the bounds proposed by Mohar \cite{Mohar91} and Lu \cite{Lu07}. From complexity point of view, the connection-graph-stability score can be interpreted as weighted edge-betweenness-centrality measure were each path is weighted by its length. By this interpretation, the connection-graph-stability scores can be calculated in polynomial time for each edge using slightly modified version of Brandes \cite{Brandes01} algorithm which has $O(NE)$ computational complexity.

This work has been supported by Swiss NSF through Grants No 200020-117975/1 and 200021-112081/1.

%In our simulations, following procedure was investigated. In initial step, a random tree-graph with $50$ nodes (and obviously with $49$ edges) was considered. Then, the potential possible edges ($1176$ edges) were added to the graph in random order (i.e. the last graph is a complete graph of order $50$). In each step, algebraic connectivity and the mentioned lower bounds were calculated. Figure (\ref{AC_LB}) presents the average of calculated algebraic
%connectivity, proposed lower bounds, and Lu and Mohar lower bounds in each step over $100$ different initial tree-graphs.
%
%\begin{figure}[tbh]
%\centering\includegraphics[width=12cm]{AC_LB.pdf}%{Lower_Bounds.jpg}%\vspace{-0.5cm}
%\caption{Algebraic connectivity and corresponding lower bounds for the graphs constructed of 100 random
%graphs of order 50. Each experiment was started with a random tree-graph and after adding the potential possible edges ended as a complete graph} \label{AC_LB}
%\end{figure}

%\section{Conclusion}
%In this paper we introduced a lower bound for algebraic connectivity
%of a graph based on maximum Connection Graph Stability number of the
%graph which is the sum of the length of paths that cross an edge. We
%proved that this new lower bound is always tighter than some other
%well-known bounds.

% \appendix

% \section{}
% \label{}

\end{document}